\documentclass[11pt,reqno]{amsart}

\usepackage{amsmath}
\usepackage{amsfonts}
\usepackage{amssymb}
\usepackage{amsthm}
\usepackage{amscd}

\usepackage{latexsym}
\usepackage{mathrsfs}
\usepackage{psfrag}
\usepackage[dvips]{epsfig}
\usepackage{epsfig}

\swapnumbers 

\theoremstyle{plain}                              % Italicized                
\newtheorem{thm}{Theorem}[section]

\newtheorem*{defn}{Definition}
\newtheorem{lem}[thm]{Lemma}
\newtheorem{cor}[thm]{Corollary}

\newtheorem*{mainthm}{Theorem}               % Or Theorem A, etc.
\newtheorem*{question}{Question}

\theoremstyle{definition}                         % Non-italicized

\newtheorem*{remark}{Remark} 

\theoremstyle{remark}                             % Non-italicized &
\newcommand{\R}{\mathbb{R}}                     % Reals
\newcommand{\ms}[1]{\mathscr{#1}}               % Script

\newcommand{\veps}{\varepsilon}

\DeclareMathOperator{\diam}{diam}

\DeclareMathOperator{\Ker}{Ker}

\DeclareMathOperator{\dR}{dist}

\providecommand{\norm}[1]{\left\lVert#1\right\rVert}       % Norm 

\providecommand{\abs}[1]{\left\lvert#1\right\rvert}        % Absolute value

\providecommand{\del}{\partial}

\begin{document}

\title[]{A lower bound on the subriemannian distance for
    H\"older distributions}

\author[S. N. Simi\'c]{Slobodan N. Simi\'c}

\address{Department of Mathematics, San
    Jos\'e State University, San Jos\'e, CA 95192-0103}

\email{simic@math.sjsu.edu}



\date{}

\begin{abstract}
  
  Whereas subriemannian geometry usually deals with smooth horizontal
  distributions, partially hyperbolic dynamical systems provide many
  examples of subriemannian geometries defined by non-smooth (namely,
  H\"older continuous) distributions. These distributions are of great
  significance for the behavior of the parent dynamical system. The
  study of H\"older subriemannian geometries could therefore offer new
  insights into both dynamics and subriemannian geometry. In this
  paper we make a small step in that direction: we prove a
  H\"older-type lower bound on the subriemannian distance for H\"older
  continuous nowhere integrable codimension one distributions. This
  bound generalizes the well-known square root bound valid in the
  smooth case.

\end{abstract}

\maketitle

%%%%%%%%%%% TEXT BODY   %%%%%%%%%

\section{Introduction}
\label{sec:intro}

The purpose of this note is to prove a lower bound in terms of the
Riemannian distance for subriemannian geometries defined by
codimension one distributions that are only H\"older continuous,
generalizing the well-known square root bound valid in the smooth
case.

A \textsf{subriemannian} or \textsf{Carnot-Carath\'eodory geometry} is
a pair $(M,H)$, where $M$ is a smooth manifold $M$ and $H$ is a
\emph{nowhere integrable} distribution (a notion we will define
shortly) endowed with a Riemannian metric $g$. Classically, both $H$
and $g$ are assumed to be smooth, e.g., of class $C^\infty$; $H$ is
usually called a \textsf{horizontal distribution} (or a polarization,
by Gromov~\cite{grom96}) and is taken to be bracket-generating (see
below for a definition). Piecewise smooth paths a.e. tangent to $H$
are called \textsf{horizontal paths}. The \textsf{subriemannian
  distance} is defined by
\begin{displaymath}
  d_H(p,q) = \inf\{ \abs{\gamma} : \gamma \ \text{is a horizontal path from $p$ to $q$} \},
\end{displaymath}
where the length $\abs{\gamma}$ of a horizontal path $\gamma : I \to
M$ is defined in the usual way, $\abs{\gamma} = \int_I
g(\dot{\gamma}(t),\dot{\gamma}(t))^{1/2} \: dt$.

Recall that a smooth distribution $H$ is called \textsf{bracket
  generating} if any local smooth frame $\{ X_1,\ldots,X_k \}$ for
$H$, together with all of its iterated Lie brackets spans the whole
tangent bundle of $M$. (In PDEs, the bracket-generating condition is
also called the H\"ormander condition.)  By the well-known theorem of
Chow and Rashevskii~\cite{mont02}, if $H$ is bracket generating, then
$d_H(p,q) < \infty$, for all $p, q \in M$.

However, $H$ does not have to be bracket generating, or even smooth,
to define a reasonable subriemannian geometry. Important examples of
such distributions come from dynamical systems. Let $f : M \to M$ be a
$C^r$ ($r \geq 1$) \textsf{partially hyperbolic diffeomorphism} of a
smooth compact Riemannian manifold $M$. This means that $f$ is like an
Anosov diffeomorphism, but in addition to hyperbolic behavior (i.e.,
exponential contraction and exponential expansion), $f$ also exhibits
weakly hyperbolic or non-hyperbolic behavior along certain tangent
directions. More precisely, $f$ is partially hyperbolic if the tangent
bundle of $M$ splits continuously and invariantly into the stable,
center, and unstable bundles (i.e., distribution),
\begin{displaymath}
  TM = E^s \oplus E^c \oplus E^u,
\end{displaymath}
such that the tangent map $Tf$ of $f$ exponentially contracts $E^s$,
exponentially expands $E^u$ and this hyperbolic action on $E^s \oplus
E^u$ dominates the action of $Tf$ on $E^c$ (see
\cite{pesin04,ps+04}). The stable and unstable bundle are always
uniquely integrable, giving rise to the stable and unstable
foliations, $W^s, W^u$. 

A partially hyperbolic diffeomorphism is (rather unfortunately) called
\textsf{accessible} if every two points of $M$ can be joined by an
$su$-path, i.e., a continuous piecewise smooth path consisting of
finitely many arcs, each lying in a single leaf of $W^s$ or a single
leaf of $W^u$. If $f$ is accessible, then the distribution $E^s \oplus
E^u$ is clearly non-integrable. Dolgopyat and
Wilkinson~\cite{dolg+wilk+03} proved that in the space of $C^r$ ($r
\geq 1$) partially hyperbolic diffeomorphisms there is a $C^1$ open
and dense set of accessible diffeomorphisms. Since invariant
distributions of a partially hyperbolic diffeomorphism are in general
only $C^\theta$, for some $\theta \in (0,1)$ (see
\cite{hps77,pesin04}), this result shows that there is an abundance of
subriemannian geometries defined by non-smooth, H\"older continuous
distributions. Since these distributions play a crucial role in the
dynamics of partially hyperbolic systems and the accessibility
property is frequently key to ergodicity (cf., e.g., ~\cite{ps+04}), a
better understanding of subriemannian geometries defined by H\"older
distributions could be very useful in dynamical systems. This paper
provides a small initial step in that direction.

\begin{figure}[tbp]
\centerline{
	\psfrag{g}[][]{$\gamma$}
        \psfrag{U}[][]{$U$}
        \psfrag{L}[][]{$L$}
        \psfrag{p}[][]{$p$}
        \psfrag{q}[][]{$q$}
        \includegraphics[width=0.4\hsize]{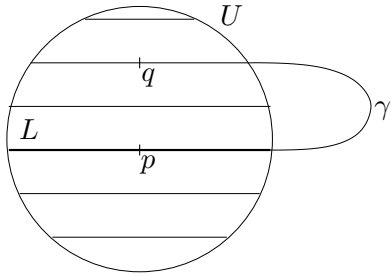}}
      \caption{Every horizontal path $\gamma$ connecting $p$ and $q$
        has to leave $U$.}
\label{fig:example}
\end{figure}

When $H$ is not smooth, the definition of nowhere integrability is
more subtle than in the smooth case. One should keep in mind the
following example (we thank an anonymous referee for bringing it up):
suppose $H$ is smooth and integrable in a neighborhood of a point $p$
and non-integrable elsewhere. Since every two points can still be
connected by a horizontal path, $d_H$ is finite everywhere. Let $L$ be
the leaf through $p$ of the local foliation of $U$ tangent to
$H$. Then unless $q$ lies in $L$, $d_H(p,q)$ is bounded away from zero
by a uniform constant, making the subriemannian distance discontinuous
with respect to the Riemannian distance for any Riemannian metric on
$M$. See Figure~\ref{fig:example}.

The following definition prohibits this type of behavior.

\begin{defn}
  A distribution $H$ on $M$ is called \textsf{nowhere integrable} if
  for every $p \in M$ and every $\veps > 0$, there exists a
  neighborhood $U$ of $p$, such that any point in $U$ can be connected
  to $p$ by a horizontal path of length $< \veps$.
\end{defn}

From now on, we will assume that \textit{the Riemannian metric on $H$
  is the restriction of a Riemannian metric from the ambient manifold}
$M$. If $H$ is a smooth bracket-generating distribution, the relation
between the subriemannian distance $d_H$ and the Riemannian distance,
which we denote by $\dR$, is well understood and is characterized by
the Ball-Box theorem~\cite{mont02, grom96}. When $H$ is of codimension
one and bracket-generating, this theorem states that in the vertical
direction, i.e., along any short smooth path $\gamma$ transverse to
$H$, $d_H$ is equivalent to $\sqrt{\dR}$. That is, there exist
constants $a, b > 0$ such that for all $p, q$ on $\gamma$,
\begin{equation}   \label{eq:rel}
  a \sqrt{\dR(p,q)} \leq d_H(p,q) \leq b \sqrt{\dR(p,q)}.
\end{equation}
In the horizontal direction, i.e., along any horizontal path, $d_H$
and $\dR$ are clearly equivalent. This means that subriemannian
geometry is non-isotropic: it behaves differently in different
directions. Subriemannian spheres are far from being ``round''. In the
Heisenberg group, for instance, subriemannian spheres
look like an apple~\cite{mont02}.

In \cite{grom96}, M. Gromov gave a short and elegant proof (without
details) of the lower bound for $d_H$, i.e., the left hand side of
\eqref{eq:rel}. His proof uses the assumption that $H$ is $C^1$ only
in the following way: if $\alpha$ is a 1-form such that $\Ker(\alpha)
= H$, then
\begin{equation}      \label{eq:stokes}
  \abs{\int_{\del D} \alpha} \leq K \abs{D},
\end{equation}
for every $C^1$ immersed 2-disk with piecewise $C^1$ boundary, where
$K$ is a constant independent of $D$. Here, $\abs{D}$ denotes the area
of $D$. If $H$ (hence $\alpha$) is $C^1$, this follows directly from
the Stokes theorem.

Gromov also remarks that without the $C^1$ assumption, the square root
estimate probably fails. We will show that this is indeed the case, in
the sense that if $H$ is only $C^\theta$, for some $0 < \theta < 1$,
and nowhere integrable, then in the vertical direction, $d_H(p,q) \geq
C \dR(p,q)^{1/(1+\theta)}$. To generalize Gromov's approach to
$C^\theta$ horizontal distributions, one needs a generalization of the
estimate \eqref{eq:stokes} to forms $\alpha$ that are only
H\"older. One such estimate was recently proved in \cite{sns+dcds+09}
and states that the integral of a $C^\theta$ $k$-form $\alpha$ over
the boundary of a sufficiently small $(k+1)$-disk $D$ is bounded by a
certain multiplicative convex combination of the $(k+1)$-volume
$\abs{D}$ of $D$ and the $k$-dimensional area $\abs{\del D}$ of its
boundary; see Theorem~\ref{thm:inequality} in the next section.

Our main result is the following:

\begin{mainthm}

  Suppose that $H$ is a nowhere integrable codimension one distribution of
  class $C^\theta$, for some $0 < \theta < 1$, on a smooth compact
  Riemannian manifold $M$. Assume that the Riemannian metric on $H$ is
  the restriction of the ambient metric from $M$. Then there exists a
  constant $\varrho > 0$ such that for any $C^1$ path $\gamma$
  transverse to $H$ and for every two points $p, q$ on $\gamma$ with
  Riemannian distance less than $\varrho$, we have
  \begin{equation}
    d_H(p,q) \geq C \dR(p,q)^{\frac{1}{1+\theta}},   \label{eq:sr}
  \end{equation}
  where $C > 0$ is a constant that depends only on $H$ and $\gamma$.

\end{mainthm}

\begin{remark}
  \begin{itemize}

  \item[(a)] Observe that if $\dR(p,q) = \veps$ is small, then
    $\veps^{1/2} \gg \veps^{1/(1+\theta)}$, which means that the lower
    bound on $d_H(p,q)$ is tighter for $C^1$ distributions than for
    $C^\theta$ ones.

  \item[(b)] It is clear that in any horizontal direction, $d_H$ is
    equivalent to $\dR$.

  \end{itemize}

\end{remark}

\section{Auxiliary results}
\label{sec:prelim}

The main tool in the proof will be the following inequality.

\begin{thm}[Theorem A, \cite{sns+dcds+09}]     \label{thm:inequality}
    
  Let $M$ be a compact manifold and let $\alpha$ be a $C^\theta$
  $k$-form on $M$, for some $0 < \theta < 1$ and $1 \leq k \leq
  n-1$. There exist constants $\sigma, K > 0$, depending only on $M$,
  $\theta$, and $k$, such that for every $C^1$-immersed $(k+1)$-disk
  $D$ in $M$ with piecewise $C^1$ boundary satisfying $\max\{
  \diam(\del D), \abs{\del D} \} < \sigma$, we have
  \begin{displaymath}
    \abs{ \int_{\del D} \alpha} \leq K \norm{\alpha}_{C^\theta} 
    \abs{\del D}^{1-\theta} \abs{D}^\theta.
  \end{displaymath}
\end{thm}

As before, $\abs{D}$ denotes the $(k+1)$-volume of $D$ and $\abs{\del
  D}$, the $k$-volume (or area) of its boundary.

\begin{remark}
  If $k = 1$, then $\diam(\del D) \leq \abs{\del D}$, so the
  assumption $\diam(\del D) < \sigma$ is superfluous.
\end{remark}

The H\"older norm of $\alpha$ on $M$ is defined in a natural way as
follows. Let $\ms{A} = \{ (U,\varphi) \}$ be a finite $C^\infty$ atlas
of $M$. We say that $\alpha$ is $C^\theta$ on $M$ if $\alpha$ is
$C^\theta$ in each chart $(U,\varphi)$; i.e., if $(\varphi^{-1})^\ast
\alpha$ is $C^\theta$, for each $(U,\varphi) \in \ms{A}$. We set
\begin{displaymath}    
  \norm{\alpha}_{C^\theta} = \max_{(U,\varphi) \in \ms{A}}
  \norm{(\varphi^{-1})^\ast \alpha}_{C^\theta(\varphi(U))}.
\end{displaymath}
For a $C^\theta$ form $\alpha = \sum a_I dx_I$ defined on an open
subset $U \subset \R^n$, we set
\begin{displaymath}
  \norm{\alpha}_{C^\theta(U)} = \max_I \norm{a_I}_{C^\theta},
\end{displaymath}
and for a bounded function $f : U \to \R$,
\begin{displaymath}
  \norm{f}_{C^\theta(U)} = \norm{f}_\infty + \sup_{x \neq y} 
  \frac{\abs{f(x) - f(y)}}{\abs{x-y}^\theta}.
\end{displaymath}

We will also need the following version of the solution to the
isoperimetric problem ``in the small''.

\begin{lem}[\cite{gromov+83}, {\bf Sublemma 3.4.B'}]  \label{lem:gromov}
  For every compact manifold $M$, there exists a small positive
  constant $\delta_M$ such that every $k$-dimensional cycle $Z$ in $M$
  of volume less than $\delta_M$ bounds a chain $Y$ in $M$, which is
  small in the following sense:
  \begin{itemize}
  \item[(i)] $\abs{Y} \leq c_M \abs{Z}^{(k+1)/k}$, for some constant
    $c_M$ depending only on $M$;
  \item[(ii)] The chain $Y$ is contained in the $\varrho$-neighborhood
    of $Z$, where $\varrho \leq c_M \abs{Z}^{1/k}$.
  \end{itemize}
\end{lem}

The following corollary is immediate.

\begin{cor}    \label{cor:gromov}
  If $\Gamma$ is a closed piecewise $C^1$ path in a compact manifold
  $M$ with $\abs{\Gamma} < \delta_M$, then there exists a $2$-disk $D
  \subset M$ such that $\del D = \Gamma$,
  \begin{equation}     \label{eq:tilde}
    \abs{D} \leq c_M \abs{\Gamma}^2,
  \end{equation}
  and $D$ is contained in the $\varrho$-neighborhood of $\Gamma$,
  where $\varrho \leq c_M \abs{\Gamma}$.

\end{cor}

\section{Proof of the theorem}
\label{sec:proof}

We follow Gromov's proof in \cite{grom96}, p. 116. Since the statement
is local and concerns arbitrary directions transverse to $H$, we can
assume without loss of generality that $\gamma$ is a unit speed
Riemannian geodesic.

Also without loss, we can assume that $H$ is transversely
orientable. If not, pass to a double cover of $M$. Let $X$ be a unit
vector field everywhere orthogonal to $H$ and define a 1-form $\alpha$
on $M$ by
  \begin{displaymath}
    \Ker(\alpha) = H, \qquad \alpha(X) = 1.
  \end{displaymath}
  Define another 1-form $\alpha_\gamma$ on $M$ by requiring that
  \begin{displaymath}
    \Ker(\alpha_\gamma) = H, \qquad \alpha_\gamma(\dot{\gamma}) = 1.
  \end{displaymath}
  Since $H$ is $C^\theta$, so are $\alpha$ and $\alpha_\gamma$. (It is
  reasonable to think of $\norm{\alpha}_{C^\theta}$ as the H\"older
  norm of $H$.) Since $\alpha$ and $\alpha_\gamma$ have the same
  kernel, there exists a function $\lambda$ such that $\alpha =
  \lambda \alpha_\gamma$. Writing $\dot{\gamma} = c X + w$, for some
  $w \in H$ and taking the Riemannian inner product with $X$, we
  obtain $\lambda = c = \cos \sphericalangle(\dot{\gamma}, X) = \sin \phi$,
  where $\phi$ is the angle between $\dot{\gamma}$ and $H$. Evaluating
  both sides of $\alpha = \lambda \alpha_\gamma$ at $\dot{\gamma}$, we
  see that $\lambda = \sin \phi$. Thus $\alpha = \sin \phi \:
  \alpha_\gamma$. Set $\phi_0 = \min \phi$. By assumption $\phi_0 >
  0$.

  Let $\tau = \min(\delta_M,\sigma)$, where $\delta_M, \sigma$ are the
  constants from Theorem~\ref{thm:inequality} and
  Lemma~\ref{lem:gromov}.

  Observe that nowhere integrability of $H$ implies that for every $p
  \in M$, the function $q \mapsto d_H(p,q)$ is continuous with respect
  to the Riemannian distance $\dR$. By compactness of $M$, it follows
  that $q \mapsto d_H(p,q)$ is in fact uniformly continuous relative
  to $\dR$. Therefore, there exists $\eta > 0$ such that for all $p, q
  \in M$, $d_H(p,q) < \tau/2$, whenever $\dR(p,q) < \eta$.

  Set
  \begin{displaymath}
    \varrho = \min \left( \eta, \frac{\tau}{2}, 
      \frac{1}{2^{\frac{1+\theta}{\theta}} c_M K^{1/\theta} 
        (\sin \phi_0)^{-1/\theta} \norm{\alpha}_{C^\theta}^{1/\theta}} \right).
  \end{displaymath}
  Observe that if $\dR(x,y) < \varrho$, then $\dR(x,y) + d_H(x,y) <
  \tau$.

  Take any two points $p, q$ on $\gamma$ satisfying $\dR(p,q) <
  \varrho$. Denote the segment of $\gamma$ starting at $p$ and ending
  at $q$ by $\gamma_0$. Let $\veps > 0$ be arbitrary but sufficiently
  small so that $\dR(p,q) + d_H(p,q) + \veps < \delta_M$. Finally, let
  $\gamma_1$ be a horizontal path from $p$ to $q$ with $\abs{\gamma_1}
  < d_H(p,q) + \veps$. (Note that a subriemannian geodesic from $p$ to
  $q$ may not exist, since $H$ is assumed to be only H\"older.)
  Define $\Gamma = \gamma_0 - \gamma_1$; $\Gamma$ is a closed
  piecewise $C^1$ path.

\begin{figure}[tbp]
\centerline{
	\psfrag{g0}[][]{$\gamma_0$}
	\psfrag{g1}[][]{$\gamma_1$}
        \psfrag{D}[][]{$D$}
        \psfrag{p}[][]{$p$}
        \psfrag{q}[][]{$q$}
        \includegraphics[width=0.4\hsize]{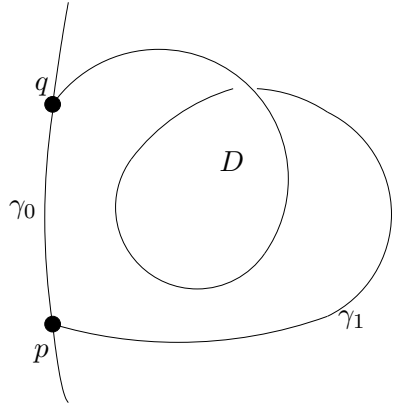}}
\caption{Vertical path $\gamma_0$ and horizontal path $\gamma_1$.}
\label{fig:paths}
\end{figure}

  Since
  \begin{displaymath}
    \abs{\Gamma} = \abs{\gamma_0} + \abs{\gamma_1} < 
   \dR(p,q) + d_H(p,q) + \veps < \tau < \delta_M,
  \end{displaymath}
  by Corollary~\ref{cor:gromov} there exists a 2-disk $D$ such that
  $\del D = \Gamma$ and $\abs{D} \leq c_M \abs{\Gamma}^2$ (see
  Fig.~\ref{fig:paths}). On the other hand, since $\tau < \sigma$, it
  follows that $\abs{\del D} = \abs{\Gamma} < \sigma$. Hence we can
  apply Theorem~\ref{thm:inequality} to $\alpha$ on $D$. Observe that
  along $\gamma_0$, we have $\alpha_\gamma = (\sin \phi)^{-1} \alpha
  \leq (\sin \phi_0)^{-1} \alpha$. Along $\gamma_1$, both $\alpha$ and
  $\alpha_\gamma$ are zero, so the inequality $\alpha_\gamma \leq
  (\sin \phi_0)^{-1} \alpha$ still holds. Therefore:
  \begin{align*}
    \abs{\gamma_0} & = \int_{\del D} \alpha_\gamma \\
    & \leq (\sin \phi_0)^{-1} \int_{\del D} \alpha \\
    & \leq K (\sin \phi_0)^{-1}\norm{\alpha}_{C^\theta} \abs{\del D}^{1-\theta} 
    \abs{D}^\theta \\
    & \leq K (\sin \phi_0)^{-1} \norm{\alpha}_{C^\theta} \abs{\del D}^{1-\theta}
    ( c_M \abs{\del D}^2)^\theta \\
    & = c_M^\theta K (\sin \phi_0)^{-1} \norm{\alpha}_{C^\theta} 
    \abs{\del D}^{1+\theta} \\
    & = c_M^\theta K (\sin \phi_0)^{-1} \norm{\alpha}_{C^\theta}
    ( \abs{\gamma_0} + \abs{\gamma_1} )^{1+\theta} \\
    & < c_M^\theta K (\sin \phi_0)^{-1} \norm{\alpha}_{C^\theta} 
    \left\{ \abs{\gamma_0} + d_H(p,q) + \veps \right\}^{1+\theta}.
  \end{align*}
  Taking the $(1+\theta)$-th root of each side and regrouping, we
  obtain
  \begin{align*}
    \abs{\gamma_0}^{\frac{1}{1+\theta}} - c_M^{\frac{\theta}{1+\theta}} 
    & K^{\frac{1}{1+\theta}}  (\sin \phi_0)^{-\frac{1}{1+\theta}} 
     \norm{\alpha}_{C^\theta}^{\frac{1}{1+\theta}} 
    \abs{\gamma_0} 
    < c_M^{\frac{\theta}{1+\theta}} K^{\frac{1}{1+\theta}} \times \\
    & \times (\sin \phi_0)^{-\frac{1}{1+\theta}} 
    \norm{\alpha}_{C^\theta}^{\frac{1}{1+\theta}} 
    \{ d_H(p,q) + \veps \}.
  \end{align*}
  Factoring $\abs{\gamma_0}^{\frac{1}{1+\theta}}$ out, the left-hand
  side becomes
  \begin{displaymath}
    \abs{\gamma_0}^{\frac{1}{1+\theta}} \left\{ 1 - c_M^{\frac{\theta}{1+\theta}} 
      K^{\frac{1}{1+\theta}} (\sin \phi_0)^{-\frac{1}{1+\theta}} 
      \norm{\alpha}_{C^\theta}^{\frac{1}{1+\theta}} 
      \abs{\gamma_0}^{\frac{\theta}{1+\theta}} \right\}.
  \end{displaymath}
  It is not hard to see that the assumption $\abs{\gamma_0} = \dR(p,q)
  < \varrho$ implies that the quantity in the curly braces is $\geq
  1/2$. Therefore,
  \begin{displaymath}
    \frac{1}{2} \abs{\gamma_0}^{\frac{1}{1+\theta}} < c_M^{\frac{\theta}{1+\theta}} 
    K^{\frac{1}{1+\theta}} (\sin \phi_0)^{-\frac{1}{1+\theta}}
    \norm{\alpha}_{C^\theta}^{\frac{1}{1+\theta}} \{ d_H(p,q) + \veps \}.
  \end{displaymath}
  Since $\veps > 0$ can be arbitrarily small and $\abs{\gamma_0} =
  \dR(p,q)$, we obtain
\begin{displaymath}
 \dR(p,q)^{\frac{1}{1+\theta}} \leq 2 c_M^{\frac{\theta}{1+\theta}} 
  K^{\frac{1}{1+\theta}} (\sin \phi_0)^{-\frac{1}{1+\theta}}
  \norm{\alpha}_{C^\theta}^{\frac{1}{1+\theta}} d_H(p,q).
\end{displaymath}
This completes the proof with
\begin{equation}    \label{eq:C}
  C = \left\{ 2 c_M^{\frac{\theta}{1+\theta}} 
    K^{\frac{1}{1+\theta}} (\sin \phi_0)^{-\frac{1}{1+\theta}}
    \norm{\alpha}_{C^\theta}^{\frac{1}{1+\theta}} 
         \right\}^{-1}. \qed
\end{equation}
\begin{remark}
  \begin{itemize}

  \item[(a)] It follows from \cite{sns+dcds+09} that $K$ stays bounded
    as $\theta \to 1-$. Therefore, by \eqref{eq:C} the constant $C =
    C(\theta)$ does not blow up as $\theta \to 1-$. This implies that
    if $H$ is $C^1$ (hence $C^\theta$, for all $0 < \theta < 1$), we
    recover the old quadratic estimate by letting $\theta \to 1-$.

  \item[(b)] Note also that as the angle between $\gamma$ and $H$ goes
    to zero (i.e., as $\gamma$ tends to a horizontal path), $C \to
    0$. This can be interpreted as follows: along horizontal paths,
    $d_H$ and $\dR$ are equivalent, so \eqref{eq:sr} cannot hold for
    any positive $C$.

  \end{itemize}
  
\end{remark}

The following question appears to be much harder than the one
investigated in this paper:

\begin{question}
  Is there an analogous \emph{upper} bound for $d_H$ with respect to
  $\dR$ when $H$ is $C^\theta$?
\end{question}

This and essentially all other questions of H\"older subriemannian
geometry remain wide open.

\bibliographystyle{amsalpha} 
                
%\bibliography{/Users/sns/Research/Articles/Biblio/master.bib}

\end{document}